\documentclass[12pt,english]{article}
\usepackage{graphicx,psfrag,amsfonts}

\newcommand{\dist}{\mbox{dist}}

\newtheorem{lemma}{Lemma}[section]
\newtheorem{theorem}[lemma]{Theorem}
\newtheorem{proposition}[lemma]{Proposition}

\newtheorem{Paragraph}[lemma]{}

\newtheorem{definition}[lemma]{Definition}
\newtheorem{remark}[lemma]{Remark}

\begin{document}

\pagestyle{myheadings} \markright {}

\footnotetext{The first author was partially supported by Project 54/001 of PDT Uruguay.}

{.}

\bigskip

{\centerline{   \Large Topological dynamics of  generic piecewise continuous }} 

{\centerline{   \Large   contractive maps in $n$ dimensions.  }}

\bigskip
 \medskip

 %% Enter the first author's name and address:
\centerline{Eleonora Catsigeras$^1$$^*$ and Ruben Budelli$^2$}
\medskip

   \centerline{$^1${\footnotesize Instituto Matem\'{a}tica, Facultad de Ingenier\'{\i}a}}
   \centerline{$^2${\footnotesize Depto. Biomatem\'{a}tica, Facultad de Ciencias}}
   \centerline {$^{1,  2}${\footnotesize Universidad de la Rep\'{u}blica, Uruguay.}}
   \medskip

 {\footnotesize  \centerline{$^*$ Corresponding author: eleonora@fing.edu.uy}}
   %\centerline{E-Addresses:  $^1$eleonora@fing.edu.uy, $^2$ruben@biomat.fcien.edu.uy}
  % {Do not forget to end the \eightpoint by the sign } here
\medskip

 %% Enter the second author's name and address:

\bigskip

    

\begin{abstract}
 We study the   topological dynamics by iterations of a piecewise continuous, non linear and locally contractive map in a real finite dimensional compact ball. We   consider those maps satisfying  the ``separation property":  different continuity pieces have disjoint images. The continuity pieces act as stable topological manifolds while the points in the discontinuity lines, separating  different continuity pieces, act as topological saddles   with an infinite expanding rate. We prove that $C^0$ generically  such systems exhibit one and at most a finite number of persistent periodic sinks  attracting all the orbits. In other words, the chaotic behaviors that this class of mappings may exhibit, are structurally  unstable and bifurcating.

\end{abstract}

\pagestyle{myheadings} \markboth { } 
{   Piecewise continuous contractive maps}

\section{Introduction}

In general, the asymptotic topological dynamics of non linear  dynamical systems with a non countably infinite set of discontinuities, is mostly unknown, particularly in large (finite) dimensions. The problem is mostly open, even if restricted  to piecewise continuous maps  acting on a compact $n$-dimensional topological manifold, and such that uniformly contract   each continuity piece.

Roughly speaking, if the map is piecewise continuous and contractive in each of its continuity pieces, then the discontinuity points, where the map can be extended but not uniquely defined, can be understood as having a   topological local splitting of saddle type, with an infinite expansion along a virtual unstable set, topologically  transversal to the stable  space where the locally contractive map $f$ acts. In fact, $f$ transforms two  points that are arbitrarily near the discontinuity, in two different images whose distance is bounded away from zero. So the virtual \lq \lq expansion rate" along that fictitious unstable set is infinite.  The continuity pieces can be translated as the  topological stable manifolds of this hyperbolic-like splitting: by hypothesis the map is locally contractive. i.e. uniformly contractive when restricted to each of its continuity pieces. Therefore,  it is hoped that the dynamics will be hyperbolic-like and thus sensible to initial conditions (topologically chaotic), if the set of discontinuities had a non empty maximal invariant subset. This  idea of a chaotic set exhibited by locally contractive piecewise $C^0$ maps, is precisely defined and proved in (\cite{Yo}).

On the other hand, if the set of discontinuities is not infinitely visited, then   the future iterates  of small neighborhoods of most orbits finally  contract forever. Therefore,  assuming a finite number of   continuity contractive pieces, and applying   the classical fixed point theorem of contractions on compact metric spaces, the limit set will be composed by a finite number of periodic attractors. In fact, we will precisely state and prove this result in Lemma \ref{lema} of this paper.

Summarizing,  the thesis that the  piecewise continuous and locally contractive systems may exhibit a chaotic attractor, or otherwise only a finite number of periodic attractors, seems natural. The main   question we pose now is:
 
 \begin{Paragraph}{  \bf Open Question.}
 \label{openquestion}
 How often the two different asymptotic dynamics, i.e. the chaotic and the periodic  behaviors respectively,  appear  in the space of piecewise continuous locally contractive maps in $n$ dimensions?
 \end{Paragraph}

Before stating the main result  in this paper that contributes to answer the question above, let us show some other motivation that lead us to study those classes of discontinuous dynamical systems.
The abstract   results about the dynamics of piecewise continuous maps in large finite dimensions may be applicable  for instance, in the theoretical understanding of the dynamics of idealized networks  composed by a large integer number  $n  $ of mutually coupled oscillators, or other attractors, with application to engineering   and secure communications (\cite{yangChua}, \cite{LuChen}). In particular the results can be used as mathematical deterministic models of very large idealized networks (for instance $n > 10^{10}$ in \cite{villa}, \cite{braun3})  of simplified biological or artificial neurons, each one behaving  as a linear or non linear oscillator,  that are mutually coupled by instantaneous inhibitory or excitatory synapsis. Those neuronal networks, if were evolving deterministically on the real time $t$, can be studied as a discrete piecewise continuous dynamical system, by means of the Poincar\'{e} first return map $f$ to a transversal section to the flow (\cite{MirolloSrogatz}, \cite{0}, \cite{Yo}).  The return map $f$  has discontinuities, due to the instantaneous synaptic jumps that couple the different oscillators. Besides, the return map $f$ is   locally contractive (\cite {0},   \cite{Yo}), if each non linear oscillator is electrically modeled by a dissipative circuit, which is  reset each time that its potential arrives to a threshold level. Numerical experiments show that, except exceptionally, the set of discontinuities (which are chaos generators) is not infinitely visited,  and the system of coupled oscillators approach asymptotically  to periodic limit cycles. In fact, for instance in \cite{braun3}, the computer simulated system modeling the dynamics of the
network of a large number of coupled neurons, shows that the global network is driven to different synchronization states that are significantly different from the original oscillatory dynamics of the individual cells. On the other hand, a bifurcating chaotic  system  of large neuronal networks (showing a disordered bursting pattern) appears in computer assisted numerical experiments, as the transition   between different periodic dynamics (ordered tonic patterns) (\cite{braun2}).

 A general theoretic   answer to the question in \ref{openquestion}, is unknown for large dimension  $n$. Nevertheless, in some particular cases, there   exist partial answers, proved classically, fitting with the computer assisted experimental results cited above. For instance,   Bruin and Deane \cite{bruin} have  proved that the contractive piecewise continuous affinities in a compact region of the plane, exhibit periodic behavior, for
Lebesgue almost every  value of a finite number of real parameters. More generally, in any dimension $n \geq 2$, C\'{e}ssac  \cite{cessac} has proved a similar theorem, for affine piecewise maps, that besides are modeling  discrete neuronal networks with synaptical instantaneous coupling on regular time intervals.

In this paper we  contribute to   answer partially the open question posed in \ref{openquestion}  including the dynamics of those maps that are not local affinities. In section \ref{seccionperiodico} we prove the following Theorem   \ref{Teorema1}:

\em Locally contractive piecewise continuous maps with the separation
property $C^0$ generically exhibit   only periodic asymptotic behavior,
 with up to a finite number of
periodic sinks that are persistent under small $C^0$ perturbations of the
map. \em

Generic  systems have a  strong $C^0$ topological meaning in this paper: they include a $C^0$ open and dense family of piecewise continuous systems. The parameters space is not reduced to a real space of finite dimension, as in \cite{bruin}, \cite{cessac}, but it is the complete functional space of all the piecewise continuous systems that are locally contractive and have the separation property. We believe that this point of view, is particularly important for the applications to mathematical dynamical models provided from other sciences. In fact, the dynamical periodic features that we prove in Theorem \ref{Teorema1} are generic and persistent, not only when moving the real parameters preserving the \em same \em model format, but also when $C^0$ perturbing \em the model itself.  \em For instance, the conclusions hold for some discrete dynamical systems obtained by iteration of the time one or of the return map to a Poincar\'{e} section, of a flow satisfying (in each of its continuity pieces)  finitely many ordinary autonomous differential equations, depending on a finite number of real parameters. But they also hold  for all the other systems obtained when   slightly changing  in all the non countably infinitely many ways, the formulae of those differential equations.

Other consequence of Theorem \ref{Teorema1}, is that the chaotic dynamics is non generic, and then it is   structurable unstable: they are destroyed if the system is perturbed, even if the perturbation is arbitrarily small. We refer to those as bifurcating systems.
 
 The thesis of Theorem \ref{Teorema1} opens the following unsolved question: is  the periodic phenomenon of locally contractive and piecewise $C^r$ maps in $n \geq 2$ dimensions,  $C^r$ dense, for $r \geq 1$? We note that our arguments to prove   the $C^0$ denseness do  not work for more regular  perturbations. 
 
 Summarizing the statement and proof of Theorem \ref{Teorema1}: \lq \lq It is essentially an extension of the classic Contraction Mapping Theorem from continuous contractive maps, to piecewise continuous   contractive (local) homeomorphisms. The weakened hypothesis however lead to a weakened conclusion, namely that the limit  set consists of a finite number of periodic orbits, rather than one fixed point, and that this behavior is only $C^0$ generic, rather than certain.  Nevertheless, this weakened conclusion is no less important for the study of dynamical systems."\footnote{Personal communication from a reader of the preliminary version of this paper.}  Besides,  the weakened conclusion is optimal, due to the fact that the bifurcating chaotic non periodic behavior  in the $C^0$ topology does exist.

    The separation property, i.e. different continuity pieces of $f$ have images whose closures are pairwise disjoint,  will play a fundamental role to prove the thesis of $C^0$- density    of the periodic behavior in the Theorem \ref{Teorema1}. But is is not necessary   for proving the persistence of the sinks in an open $C^0$ neighborhood of each map $f$ exhibiting a periodic behavior. We note that the hypothesis of separation is also unnecessary, if substituted by  other  assumptions, to obtain certain prevalence of the periodic behavior.  For instance, if $f$ is an affinity in each of its continuity pieces, then the separation property is not needed  \cite{bruin, cessac}. We conjecture that the separation property is also unnecessary to obtain a $C^2$ genericity of the periodic behavior, if the system is not piecewise-affine but is piecewise $C^2$ and thus exhibits, in each of its continuity pieces, a property of bounded distortion of the derivative of the $k-$th iterate, uniformly for all $k \geq 1$. (Note that affine maps, in particular, have zero distortion in  the derivative of $f^k$).

\section{Definitions.} \label{definicionesabstractas}

\begin{definition}
\label{particion}

\em
 Let $B \subset {\mathbb{R}^{n}}$  be a compact set,  homeomorphic  to a compact ball. In particular $B $ is connected.
A \em finite partition \em of $B$ is a finite collection
$\{B_i\}_{1\leq i \leq m}$  of
compact non empty sets $B_i \subset B$,
 such
that
 $ \bigcup _{1\leq i \leq m}B_i = B$  and  $\mbox{ int }B_i \, \cap
 \mbox{ int }B_j \, = \, \emptyset $, for $ i \neq j $.
Denote $S =  \bigcup _{i \neq j } B_i \cap B_j = \bigcup_{i=1}^m \partial B_i \,$ and call $S$
the \em separation line, \em  or \em line of discontinuities, \em (although it is not a line in the usual sense, but the union
of the  topological frontiers of $B_i$).

We endow the set $B \subset {\mathbb{R}^n}$ with \em any \em metric \lq\lq$\mbox{dist}$" (which is nor necessarily derived from an inner product neither from a norm in $\mathbb{R}^n$), but such that it still induces  the   usual  topology of $B  $ as a subset of the Hilbert space $\mathbb{R}^n$. The  freedom of the results in this section (particularly of the thesis of Theorem \ref{Teorema1}), from the chosen metric $\mbox{dist}$ in $B$, is relevant for the applications to some physical piecewise continuous dynamical systems. In fact, for instance in \cite{Yo}, it is proved that the mathematical model of a large network of inhibitory neurons, evolving deterministically in time $t \in \mathbb{R}^+$, can be described through its first return map $f$ to a Poincar\'{e} section, and this $f $ is a  piecewise continuous dynamical system, locally contractive respect to certain metric \lq \lq$\mbox{dist}$". This metric is usually constructed as being adapted to the dynamics, as for hyperbolic sets, and so, unless the system were linear, it is not the derived distance, nor from an inner product neither from a norm in $\mathbb{R}^n$. So \lq \lq$\mbox{dist}$" is not in general the usual distance induced from $\mathbb{R}^n$) on $B$. The only condition, obtained in \cite{Yo} as a thesis, and used in the topological assertions all along the proof of the Theorem \ref{Teorema1} of this paper, is that the metric \lq\lq$\mbox{dist}$" induces  the usual topology in $B$ as a subset of $\mathbb{R}^n$.
\end{definition}

\begin{definition} \em \label{piecewisecontinuous}

Given a finite partition $\{B_i\}_{1\leq i \leq m} $ of $B$,
we call $F$ a   \em  piecewise
 continuous map on $(B, {\cal P})$ with the separation property \em
  if $F$ is a finite family  $
 F = \{f_i\}_{1\leq i \leq m}$
 of  \em homeomorphisms \em $f_i: B_i \mapsto f_i(B_i) \subset int(B) $, such that $f_i (B_i) \cap
 f_j (B_j) = \emptyset$ if $i \neq j$. We note that $F$ is multi-defined in the separation line $S$.

  Each $B_i$ shall be called a \em continuity piece \em of $F$.

  \end{definition}

   \begin{remark}
   \label{remarkFinversa}
    \em  A piecewise continuous map $F$ with the separation property is globally one to one because it is an homeomorphism in each continuity piece and two different continuity pieces have disjoint images. Therefore $F^{-1}$ exists, uniquely defined in each point of  $F(B) = \bigcup_i f_i(B_i)$. In fact:

   For any point $x \in \bigcup_i f_i(B_i)$, its backward first iterate is uniquely defined as  $F^{-1}(x) = f_i^{-1}(x)$, where $i$ is the unique index value such that $x \in f_i(B_i)$.

   Nevertheless $F^{-1}$ is  not necessarily one to one because  $F$ is multidefined in $S = \bigcup_{i \neq j} (B_i \cap B_j)$.

    $F^{-1}$ is continuous in $F(B)$, because $F^{-1}|_{f_i(B_i)} = f_i^{-1}$ and $f_i$ is an homeomorphism due to the Definition \ref{piecewisecontinuous}.
   \end{remark}

\begin{definition} \em  \label{definicioncontraccionlocal}

We say that $F$ is  uniformly \em locally contractive  \em if there exists a constant $0 < \lambda < 1$, called an \em uniform contraction rate for $F$, \em  and a metric $\dist$ in $B$, such that
$\dist (f_i(x),f_i(y)) \leq \lambda \dist (x,y)$, for all $x$ and $y$ in the same $B_i$,
for all $1 \leq i \leq m $ .

 \end{definition}

Given a point $x \in B$,     its image set  is $F(x) = \{f_i(x):   x \in B_i \}$.
If $H \subset B$, its image set is $F(H) = \bigcup _{x \in H} F(x)$.
We have that $B \supset F(B) \supset \ldots F^k(B) \supset \ldots $.

The second iterate of the point $x \in B$ is the set $F^2(x) = F (F(x))$. In general, for all $j \geq 1$ we define  the $j-$th. iterate of $x \in B$ as the set $F^{j}(x) = F^{j-1}(F(x))$. We agree to define $F^0(x) = \{x \}$ and $F^0(H) = H$.

\begin{definition} \em \label{definicionAtomo}
For any natural number $k \geq 1$, we define the \em atoms of generation $k$ \em as the sets
 $$f_{i_k}\circ \ldots \circ f_{i_2} \circ
 f_{i_1} (B_{\mathbb{I}}) $$ where $\mathbb{I} = (i_1, i_2, \ldots, i_k) \in \{1,2,\ldots, m\}^k$
 and $B_{\mathbb{I}}$ is  the subset of $B_{i_1}$ where the  composed function above is defined. (If $B_{\mathbb{I}}$ were an empty set, then the atom is empty.)  Abusing of the notation we write the atom as:
 $$f_{i_k}\circ \ldots \circ f_{i_2} \circ
 f_{i_1} (B_{i_1})$$
 \end{definition}

 We note that each atom of generation $k$ is a compact,
  not necessarily connected set, whose
 diameter is smaller than $\lambda ^k \mbox{diam} B$.

 The set $F^k(B) $ is a compact set, formed by the union of all the
 non empty atoms
of generation $k$.
There are at most $m^k$ and at least $m$ non empty atoms of generation $k$, where $m$ is the number of continuity pieces of $F$.

\begin{definition} \label{limitset}
\em

Given $x_0 \in B$, a future orbit $o^+(x_0)$ is a sequence of points $\{x_i\}_{i \geq 0}$, starting in $x_0$, such that $x_{i+1} \in F(x_i) \; \; \forall \, i \geq 0$. Due to the multi-definition of $F$ in the separation line $S$, the points of $S$ and those that eventually fall in $S$ may have more than one future orbit.

A point $y$ is in the limit set $L^+(o^+(x_0))$ of a future orbit of $x_0$ if there exists $k_j \rightarrow +\infty$ such that $x_{k_j} \rightarrow y$.

The limit set $L^+(x_0)$ is the union of the limit sets of all the future orbits of $x_0$.

  The \em limit set \em $L^+(B)$ of the map $F$, also denoted as $L^+(F)$, is
 the union of the limit sets $L^+(x)$ of all the points $x \in B$. \em
 \end{definition}

 \begin{remark} \em
 \label{remarkConjLimiteinvariante} Due to the compactness of the space $B$ the limit set
 $L^+(o^+(x_0))$ of any future orbit, \em is non empty. \em It is standard to prove that \em $L^+(o^+(x_0))$ is compact \em (because it is closed in the compact space $B$). Nevertheless $L^+(x_0)$ may be not compact, if the point $x_0$ has infinitely many different future orbits.
Finally, we assert that \em $L^+(o^+(x_0))$ is invariant: \em $F^{-1} (\; L^+(o^+(x_0)) \;  ) = L^+(o^+(x_0))$.

 \end{remark}

 {\em Proof:}  Consider $y \in L^+(o^+(x_0))$. We have $y = \lim _{j \rightarrow + \infty} x_{k_j} \in F(B)$ if $k_j \geq 1$.

  $F^{-1}: F(B) \rightarrow B$ is a continuous uniquely defined function  (see Remark \ref{remarkFinversa}).  Then $x_{k_j -1} = F^{-1} (x_{k_j}) \rightarrow F^{-1}(y) $, so $F^{-1} (y) \in L^+(o^+(x_0))$ proving that $$F^{-1} (\; L^+(o^+(x_0)) \;  ) \subset L^+(o^+(x_0)) $$

 Let us prove the converse set inclusion: $F^{-1} (\; L^+(o^+(x_0)) \;  ) \supset L^+(o^+(x_0)) $.

 $F = \{f_i: B_i \mapsto B\}$ is defined and continuous in each of its finite number of pieces $B_i$, that are compact sets that cover $B$. Then there exists some $i \in\{1, 2, \ldots, n\}$ and a subsequence (that we still call $k_j$), such that $$y = \lim_{j \rightarrow + \infty} x_{k_j} \in B_i,  \; \; \; \forall j \geq 0: \; \; x_{k_j} \in B_i, \; \; x_{k_j + 1} = f_i (x_{k_j}), $$ $$ f_i(y) = \lim f_i(x_{k_j}) = \lim x_{k_j + 1} $$
 We conclude that there exists $y_1 = f_i(y) \in F(y)$ such that $y_1 \in L^+(o^+(x_0))$. In other words, $y \in F^{-1}(L^+(o^+(x_0))$. This last assertion was proved for any $y \in L^+(o^+(x_0))$. Therefore $ L^+(o^+(x_0)) \subset F^{-1} (\; L^+(o^+(x_0)) \;  ) $ as wanted. $ \; \Box$
 \begin{definition}

\em
 We say that a point $x$ is \em  periodic of period $p$  \em
 if there exists a first natural number $p\geq 1$   such that $x \in F^{p}(x) $. This is
 equivalent to $x$ be a periodic point in the usual sense, for
  the uniquely defined map $F^{-1}$, i.e.
  $F^{-p}(x)= x$ for some first natural number $p \geq 1$.

  We call the backward orbit
 of $x$ (i.e. $\{ F^{-j}(x), j = 1,\ldots, p\}$), a periodic orbit
  with period $p$.

\end{definition}

   It is not difficult to show that   the limit set $L^+(F)$ is contained in the compact, totally
   disconnected set $K_0 = \bigcap _{k\geq1} F^k(B)$. It could be
   a Cantor set. But generically
    $K_0$  shall be the union of a finite number of periodic orbits,
    as  we shall prove in Theorem \ref{Teorema1}.

\begin{definition} \em  \label{finallyperiodic}

  We say that $F$ is \em finally periodic \em with period $p$ if the
  limit set $L^+ (F)$ is the union of \em only \em a finite number of periodic orbits
  with least common multiple of their periods equal to $p$.
 In this case we call \em limit cycles \em to the periodic orbits of  $F$.

 We call basin of attraction of each limit
 cycle  $L$ to the set of points $x \in B$ whose  limit set $L^+(x)$ is $L$.

\end{definition}

{\bf Topology in the space of piecewise continuous locally contractive maps.}

  Let ${\cal P} = \{B_i\}_{1\leq i \leq m}$ and  ${\widehat{{\cal P}}} =
  \{A_i\}_{1\leq i \leq m}$  be   finite partitions (see Definition \ref{particion})
  of the compact region $B$ with the same number $m$ of pieces.

We define the distance between ${\cal P}$ and ${\widehat{{\cal P}}}$ as
\begin{equation}
\label{definicionHausdorfDistance}
d({\cal P}, {\widehat{{\cal P}}}) = \max_{1 \leq i \leq m}\; \; \;  \mbox{Hdist} (A_i, B_i) \end{equation} where $\mbox{Hdist} (A_i, B_i)$ denotes the Hausdorff distance between the two compact sets $A_i$ and $B_i$. i.e.
$$\mbox{Hdist} (A_i, B_i) = \max \{\dist (x, B_i), \dist (y,A_i):  x \in A_i, y \in B_i  \}$$
and $\dist (x, B_i) = \min \{ \dist (x,y): y \in B_i\}$

Although it is standard to check the following properties of the distance between two partitions ${\cal P}$ and ${\widehat{{\cal P}}}$, we include their proofs for a seek of completeness:

\begin{proposition}
\label{remarkagregado}.    If $d({\cal P}, {\widehat{{\cal P}}}) < \epsilon $ then:
\em

 -$\mbox{ Hdist }(S, \widehat S) < \epsilon$, \em where $S = \cup_{i } (\partial B_i ) $ is the separation line of the partition ${\cal P} = \{B_i: 1 \leq i \leq m  \}$, and $\widehat S = \cup _{i} (\partial A_i)$ is the separation line of the partition ${\widehat{{\cal P}}} = \{A_i: 1 \leq i \leq m\}$.

- For all $i \neq j$ such that $B_i \cap A_j \neq \emptyset$, and for all $p \in B_i \cap A_j$: \em $$\dist (p, S) < \epsilon, \; \; \dist (p, \widehat S) < \epsilon.$$

\end{proposition}

{\em Proof: }  In the following proof we will use that $B$ is homeomorphic to a compact ball in $\mathbb{R}^{n}$: it is a compact and connected metric space and so, all the subsets $M \subset B$ have the following property:
$$y \not \in M \; \; \Rightarrow \; \; \; \dist (y, M) = \dist (y, \partial M)$$
where $\partial M$ is the topological frontier of $M$ as a subset of the topological space $B$.

To deduce that $\mbox{ Hdist }(S, \widehat S) < \epsilon$, recall that $S= \cup_{i=1}^m \partial B_i, \; \; \widehat S = \cup_{i=1}^m \partial A_i$. Then:
 $$\mbox{ Hdist} (S, \widehat S) = \max \{\dist (x, S), \dist (y,\widehat S):  x \in \cup_{i=1}^m \partial A_i, y \in \cup_{i=1}^m \partial  B_i         \} = $$
 $$= \max \{\min_{1 \leq j \leq m} \{\dist (x, \partial B_j)\}, \min_{1 \leq j \leq m}\{\dist (y, \partial A_j)\}:  x \in \cup_{i=1}^m \partial A_i, y \in \cup_{i=1}^m \partial  B_i         \} \leq$$
 $$\leq  \max_{1  \leq  i \leq m} \{ \max \{\dist (x, \partial B_i), \dist (y,\partial A_i):  x \in  \partial A_i, y \in \partial  B_i         \}\} =$$
 $$= \max_{1  \leq  i \leq m} \mbox{Hdist} (\partial B_i, \partial A_i)$$
 
  Therefore, to prove that $\mbox{ Hdist }(S, \widehat S) < \epsilon$, it is enough to prove that $\mbox{ Hdist} (\partial B_i, \partial A_i) < \epsilon$ for all $i \in \{1, \ldots, m\}$. 
  
  If $\partial B_i = \partial A_i$ then their Hausdorff distance is zero and thus, smaller than $\epsilon$. On the other case, there exists $p \in (\partial B_i \setminus \partial A_i ) \cup ( \partial A_i \setminus \partial B_i)$. First suppose $p \in \partial B_i, \; p \not \in \partial A_i$.
$$d({\cal P}, {\widehat{{\cal P}}}) < \epsilon \; \Rightarrow \; \; \mbox {dist} (p, A_i) < \epsilon \; \; \forall \, p \in B_i, \mbox{ in particular } \; \forall \, p \in \partial B_i  $$
If $p \not \in \partial A_i$ then $\mbox{dist}(p, A_i) = \dist (p, \partial A_i)$
$$ \Rightarrow \mbox{ dist } (p , \partial A_i) < \epsilon \; \; \; \forall \, p \in \partial B_i$$
Changing the roles of $A_i$ and $B_i$, the same argument works for
 $q \in \partial A_i \setminus \partial B_i$. So we deduce
$$\mbox{ Hdist } (\partial A_i, \partial B_i) = $$ $$= \max \{\dist (p, \partial A_i), \dist (q, \partial B_i), \; p \in \partial B_i, \; q \in \partial A_i   \} < \epsilon$$
Let us prove now the second assertion in this remark. We will only prove that $\dist (p, \widehat S) < \epsilon \; \; \forall \, p \in B_i \cap A_j$. The  inequality
$\dist (p,  S) < \epsilon$ follows from this one, changing the roles of the partitions ${\cal P}$ and ${\widehat{{\cal P}}}$.

If $p \in  B_i  \cap A_j$ then, being $i \neq j$, we deduce $$\mbox{int}A_i \cap \mbox{int} A_j = \emptyset \; \Rightarrow \; \mbox{int}A_i \cap  A_j = \emptyset \; \; \Rightarrow \; p \not \in \mbox{ int}{ A_i} \; \; \Rightarrow \; \; p \in B_i \setminus (\mbox{int} A_i)$$
$$\Rightarrow \; \; \dist (p, \mbox{int} A_i) = \dist (p , \partial A_i) \leq \mbox{Hdist} (B_i, A_i) < \epsilon $$
But $\partial A_i \subset \widehat S$, then $ \dist (p, \widehat S) \leq \dist (p, \partial A_i) $.
So we conclude $\dist (p,  \widehat S) < \epsilon$ as wanted. $ \; \; \Box$

\begin{definition} \em  \label{topologia}
  Let
   $
  F  = \{f_i:B_i \mapsto B\}_{1\leq i \leq m}$
  and $G   =\{g_i:A_i \mapsto B\}_{1\leq i \leq m}$  be    locally contractive
   piecewise continuous maps on $(B, {\cal P})$ and  $(B, {\widehat{{\cal P}}})$ respectively.
    Given
   $\epsilon >0 $
  we say that   $G $ is a
  \em $\epsilon$-perturbation
 of $F$  \em if
$$ \max_{1 \leq i \leq m} \left \| \left. (g_i - f_i) \right | _ {\displaystyle {B_i \cap A_i}} \right \|_{\displaystyle {{\cal C\;} ^0}} < \epsilon, \; \; |{ \lambda_F - \lambda_G} | < \epsilon \; \; \; \mbox { and
} \; \; d({\cal P}, {\widehat{{\cal P}}} ) < \epsilon $$
where $\lambda _F$ denotes \em the uniform contraction rate \em of $F$ in  its continuity pieces, defined in \ref{definicioncontraccionlocal}, and $\|\cdot \|_{{\cal C}^0}$ denotes the  ${\cal C}^0$ distance in the functional space of continuous functions defined in a \em compact \em set $K$: $$\| (g - f)|_K\|_{{\cal C}^0} = \max_{ x \in K} \dist (g(x), f(x))$$
\end{definition}

\begin{definition} \label{definicionPersistencia} \em
 We say that the limit cycles  of a finally periodic map $F$  (see Definition \ref{finallyperiodic})
  are \em persistent  \em if:

  For all $\epsilon^*>0$ there exists $\epsilon >0$
   such that all $\epsilon$-perturbations $G$
of $F$  are finally periodic  with the same
finite number of limit cycles (periodic orbits) than $F$, and such that each limit  cycle $L_{G}$ of $G$ has the same period and is $\epsilon^*$-near of some limit cycle $L_{F}$ of $F$ (i.e. the Hausdorff distance between $L_G$ and  $L_F$ verifies
$\mbox{Hdist}(L_G, L_F) < \epsilon^*$).

\end{definition}

\begin{definition} \em  \label{generico}

Denote ${\cal S}$ to the \em space of all the systems that are  piecewise continuous with the separation property and locally contractive, \em according with the Definitions \ref{piecewisecontinuous} and \ref{definicioncontraccionlocal}.

We say that a property $\mathbb{P}$  of the systems in ${\cal S}$ (for instance being finally periodic as we will  show in Theorem \ref{Teorema1}) is \em (topologically) generic \em
if $\mathbb{P}$ is verified, at least, by an \em open and dense \em subfamily of systems in the functional space ${\cal S}$, with the topology  in ${\cal S}$ defined in \ref{topologia}.

Precisely, being \em generic \em means:

1) The \em openness \em  condition: For each piecewise continuous map $F$ that verifies the property $\mathbb{P}$ there exist $\epsilon >0$ such that
all $\epsilon$-perturbation of $F$ also verifies $\mathbb{P}$.

2) The \em denseness \em condition:  For each piecewise continuous map $F$ that does not verify the property $\mathbb{P}$, given $\epsilon >0$, arbitrarily small, there exist some $\epsilon$-perturbation $G$ of $F$ such that $G$ verifies the property $\mathbb{P}$.

The openness  condition implies that the property $\mathbb{P}$ shall be robust under small perturbations of the system. It is robust under small changes, not only of a finite number of real parameters, but also of the functional parameter that defines the model itself. So   the system should be structurally stable. When this robustness holds, the property $\mathbb{P}$ is still observed  when the   system,  the model itself, does not stay  exactly fixed, but is changed, even  in some unknown fashion, remaining near the original one.

The density condition combined with the openness condition,  means that the only behavior that have chance to be observed under not exact experiments are those that verify the property  $\mathbb{P}$. In fact, if the system did not exhibit the property $\mathbb{P}$, then some arbitrarily small change of it, would lead it to exhibit $\mathbb{P}$ robustly. \em

\end{definition}

The denseness condition implies that if the property $\mathbb{P}$ were generic, then the opposite property   (Non-$\mathbb{P}$) has null interior in the space of ${\cal S}$ of systems, i.e. Non-$\mathbb{P}$ is not robust: some arbitrarily small change in the  system will lead it to exhibit $\mathbb{P}$. That is why we define the following:

\begin{definition} \em  \label{bifurcating}
If the property $\mathbb{P}$ is generic, we say that any system that does not exhibit $\mathbb{P}$ is \em bifurcating, \em  and  Non-$\mathbb{P}$ is a \em not persistent \em property.
\end{definition}

\section{The generic persistent periodic behavior.} \label{seccionperiodico}

\begin{theorem} \label{Teorema1}
Let $F$ be a locally contractive  piecewise continuous
map with the separation property.  Then generically $F$ is finally periodic with persistent limit cycles. \em

\end{theorem}

To prove Theorem \ref{Teorema1} we shall use the following lemmas \ref{lema} and \ref{lemma2}:

 \begin{lemma} \label{lema}
If there exists an integer $k_0 \geq 1$ such that the compact set
 $K_0 =  F ^{k_0} (B)$ does not intersect the separation line $S$ of the partition into the continuity pieces of $F$, then
 $F$ is finally periodic.

\end{lemma}

{\em Proof: }
 By hypothesis $\dist (K_0, S) = d >0$,  because $K_0$ and $S$ are disjoint compact sets. On the other hand $$K_0= F^{k_0}(B) = \bigcup _{A \in {\cal A}_{k_0}} A$$ where ${\cal A}_{k}$ for any fixed $k \geq 1 $, denotes the family of all the atoms of generation $k$ defined in \ref{definicionAtomo}.

 The  diameter $\mbox{ diam}(A)$ of
 each atom $A$ of the finite family ${\cal A}_k$, is smaller than $\mbox{diam} (B) \, \lambda^k $. Therefore it converges to zero when $k \rightarrow +\infty$. Thus, for all  $k $  large enough:

 $$\mbox{ diam} (A) \leq \frac{d}{2} \; \; \; \; \forall \; A \in {\cal A}_k$$

 It is not restrictive to suppose $k \geq k_0$. Then $A \subset F^k(B) \subset F^{k_0}(B) = K_0 \; \; \;
 \forall \; A \in {\cal A}_k. $

 We assert that each atom $A \in {\cal A}_k$ for such  $k$, is contained in the interior of some continuity piece $B_i$.
To prove this last assertion we give the following argument (P), that will  be useful also in the proof of Lemma \ref{lemma2}:

 {\bf (P)} Fix a point $x \in A$. As the continuity pieces cover the space $B$, there exists some (a priori not necessarily unique) index $i$ such that $x \in B_i$. It is enough to prove that $y  \in \    int (B_i)$ for all $y \in A$ (including $x$ itself).

 We argue  in the compact and connected metric space $B$, using the following known properties of the  metric space $B$ with the topology induced by its inclusion in $\mathbb{R}^{n}$, as a subset homeomorphic to a compact ball.

   - The triangle property.

   - The distance $\mbox{dist} (y,M)$ of a point $y \not \in M$, to a set $M \subset B$, is the same that the distance of $y$ to the topological frontier $\partial M$ of $M$  as a subset of $B$.

   -  $\mbox {dist}(y,M_1) \geq \mbox{ dist}(y,M)  $ if $M_1 \subset M$.

  We denote $B_i^c$ to the complement of $B_i$ in $B$, and in the topology relative to $B$  we denote: $\overline {(B_i^c)}$ to the closure of $B_i^c$, i.e the complement of $int(B_i)$, and $\partial B_i$ to the frontier of $B_i$ in $B$, $\partial B_i \subset S$:
 $$\dist (x,y) \leq diam (A) < d/2, $$ $$\dist (x, \overline {(B_i^c)}) = \dist (x, \partial B_i) \geq \dist (x, S) \geq d$$ $$ \dist (y, \overline {(B_i^c)}) \geq \dist (x, \overline{(B_i^c)}) - \dist (x,y)  \geq d - d/2 = d/2 >0 $$
 Therefore $y \not \in \overline{(B_i^c)}$ proving that $y \in \ int (B_i) $ as wanted. $\; \Box $ {\bf (P)}

We deduce that given an atom $A \in {\cal A}_k$, there exists a unique  natural number $i_0$ such that $A \in int (B_{i_0})$. Therefore $F(A)$ is a single atom of generation $k+1$.

  From the definition of atom in \ref{definicionAtomo}, we obtain  that any atom of generation larger than $k$ is contained in an atom of generation $k$. But each atom of generation $k$ is in the interior of a piece of continuity of the partition $\{B_i\}$.
We deduce that there exists a sequence of natural numbers $\{i_h\} _{h \geq 0}$, called \em the itinerary of  the atom $A$, \em such that \begin{equation}
  \label{equationImagenesDeAtomo}
  A \subset \ int (B_{i_0 }), \; \; F(A) = f_{i_0}(A)  \subset \ int (B_{i_1}), $$ $$ F^2(A) =  f_{i_1} \circ f_{i_0} (A) \subset \ int (B_{i_2}), \ldots\end{equation} and the successive images of the atom $A$ of generation $k$, are  single
   atoms of generation $k+1, k+2, \ldots, k+h, \ldots $.  Therefore, the successive images of the atom $A$, in the sequence (\ref{equationImagenesDeAtomo}), are contained in a sequence of atoms: $A= A_0, A_1, A_2, \ldots, A_h, \ldots, $ all of generation $k$.

   The same property holds for any of these atoms of generation $k$, and each of them is contained in the interior of a continuity piece of $F$, so $F$ is uniquely defined there and we have:
   \begin{equation}
   \label{ecuacionChain2}
   A = A_0 \subset \ int (B_{i_0}), \; \; F(A_0) \subset A_1 \subset \ int (B_{i_1}) , $$ $$F^2(A_0) \subset F(A_1) \subset A_2  \subset \ int (B_{i_2}) , \ldots, \end{equation}

   For fixed  $k$, the family of atoms of generation $k$ is finite, so we conclude that there exists two first natural numbers $0 \leq h < h+p$ such that $F^p (A_h) \subset A_h$.

   Note that,  $F^p (A_h) $ is uniquely defined as $ f_{i_{h+p}} \circ f_{i_{h+p -1}} \circ \ldots \circ f_{i_{h}}$, because we are considering sets contained in the interior of the continuity pieces of $F$.

    Due to the uniform contractiveness of $f_i$ in each of its continuity pieces, $F^p : A_h \mapsto A_h$, is uniformly contractive. The Banach Fixed Point Theorem (i.e. the Contraction Mapping Theorem), states that in a complete metric space, any uniformly contractive map from a compact set to itself, has an unique fixed point, and  all the orbits in the set converge to this fixed point in the future. Therefore, there exists in $A_h$ a periodic point $p_0$ by $F$, of period $p \geq 1$, and all the orbits with initial states in $A_h$ have the periodic orbit $L$ of $p_0$, as their  limit set.

    By construction $A_h$ contains the image of $A$ by an iterate $F^h$, uniquely defined. So we conclude that the limit set of all the points in the atom $A$ is $L$.

    The construction above can be done starting with any initial atom $A \in {\cal A}_k$. Besides, ${\cal A}_k$ is a finite family. We conclude that there exists one, and at most a finite number of periodic limit cycles, attracting all the orbits of $\bigcup_{A \in {\cal A}_k} A = F^k(B)$.

    The last assertion implies that the limit set of $B$ is formed by that finite family of periodic limit cycles, ending the proof of this lemma.
    $ \; \; \; \; \Box$

\begin{lemma}
\label{lemma2}
In the hypothesis of Lemma  \em \ref{lema}, \em the limit cycles of $F$ are persistent.
\end{lemma}

{\em Proof:}
    We shall prove that the limit cycles are persistent according to the definition \ref{definicionPersistencia}.

    The condition of the hypothesis of  Lemma  \ref{lema} is open in the topology defined in \ref{topologia}, because $K_0$ and $S$ are compact and at positive distance. Therefore, there exists $\epsilon_0>0$ such that, for all $0<\epsilon < \epsilon _0$, all $\epsilon-$perturbation $G$ of $F$, is finally periodic.

    {\bf (Q) }  We claim that, \em given ${k_0} \geq 1$ fixed such that $\dist (A, S ) \geq d >0$ for all $A \in {\cal A}_{k_0}(F)$, then there exists  $0<\epsilon< \epsilon_0$  small enough such that if $G$  is a $\epsilon-$perturbation of $F$, then there is a bijection $\Psi$ between the families ${\cal A}_{k}(F)$ and ${\cal A}_{k}(G)$, of the atoms of all generation ${k \geq 1}$ of $F$ and $G$ respectively, and besides, for some $k$  large enough, the itinerary of each of the atoms $A \in {\cal A}_{k}(F)$ is the same than the itinerary of the respective  atom $\Psi (A) = \widehat A \in {\cal A}_{k}(G)$. \em

    In fact, due to the definition of $\epsilon_0$- perturbation of $F$, the continuity pieces $B_i = B_i(F) \subset B$ and $\widehat B_i = B_i(G) \subset B$, of $F$ and $G$ respectively, are correspondent by a bijection, such that the Hausdorff distance  $$\mbox{Hdist} (B_i(F), B_i(G)) < \epsilon_0.$$

     On the other hand, for all $k \geq 1$, the atoms $A \in {\cal A}_{k}(F) $ and $\widehat A \in {\cal A}_{k}(G)$, due to the definition of atom in \ref{definicionAtomo}, are:   $$A= F^{k}(B_{\mathbb{I}}) , \; \;  \widehat A = G^{k}(B_{\widehat{\mathbb {I}}})$$  identified by  words $$\mathbb{I} = (i_1, i_2, \ldots, i_{k}) , \; \;  \widehat {\mathbb{I}} = ( \widehat i_1, \widehat i_2, \ldots, \widehat i_{k})\; \;  \in \{1, 2, \ldots, m\}^{{k}}$$

      We define the correspondence $$\Psi (A) = \widehat A \; \; \; \mbox{ if and only if } \; \; \;  \widehat {\mathbb I} = \mathbb{I}$$

      With $k= {k_0}$ fixed, we have \begin{equation}
      \label{equation111}
      \mbox{ dist } (A, S ) \geq d >0 \; \; \; \forall \; A \in {\cal A}_{k_0}(F).\end{equation} On the other hand, due to the definition \ref{topologia} of $\epsilon_1-$perturbation $G$ of $F$, we have $\|f_{i} - g_{i}\|_{C^0} <  \epsilon_1, \; \; \; \mbox{H dist} (B_i, \widehat B_i) < \epsilon_1 \; \; \; \forall \; i = 1, 2, \ldots, m$. But the finite composition of   contractive homeomorphisms depends continuously of the homeomorphisms, in the topology defined in \ref{topologia}. Then,  for ${k_0}$-fixed,  and using the definition of continuity, for given $\epsilon = d/3$ (we choose $\epsilon = d/3$), there exists $0 <\delta = \epsilon_1< \epsilon _0 $ (it is just a matter of notation, to call $\epsilon_1$ to $\delta$, for convenience in further uses), such that
      $$\|f_{i} - g_{i}\|_{C^0} < \delta = \epsilon_1, \; \; \; \mbox{Hdist} (B_i, \widehat B_i) < \delta = \epsilon_1 \; \; \; \forall \; i = 1, 2, \ldots, m \;\; \Rightarrow \; \; \; $$ $$\|f_{i_{k_0}} \circ f_{i_{{k_0}-1}} \circ \ldots \circ f_{i_{1}}(B_{i_{1}}) \; - \; g_{i_{k_0}} \circ g_{i_{k_0-1}} \circ \ldots \circ g_{i_{1}}(\widehat B_{i_{1}})\| < \epsilon = \frac{d}{3} $$ $$ \forall \; \mathbb{I} = (i_1, i_2, \ldots, i_{k_0}) \in \{1,2, \ldots, m\}^{k_0}$$
      In other words, the last statement can be reformulated as:
      $$G \mbox{ is a } \epsilon_1- \mbox{perturbation of } F \; \; \; \Rightarrow \; \; \;$$ $$  \mbox{ Hdist }(A, \Psi (A) ) < \frac{d}{3} \; \; \; \forall \; A \in {\cal A}_{k_0}(F) \; \mbox{ and } \,  \widehat A = \Psi (A) \in {\cal A}_{k_0}(G).$$
      Besides, if $\epsilon_1>0$ is chosen smaller than $d/3$, from the definition \ref{topologia} we obtain $$\mbox{ H dist} (S , \widehat S) < \epsilon_1< \frac{d}{3}$$
      where $\widehat S$ is  line of discontinuities of the piecewise continuous map $G$, which is a $\epsilon-$ perturbation of $F$.
      Joining the last two inequalities with (\ref{equation111}) and applying the triangular inequality, we deduce:
      $$\mbox{ dist } (\widehat A, \widehat S) \geq d- \frac{d}{3} - \frac{d}{3} = \frac{d}{3}  \; \; \; \forall \; A \in {\cal A}_{k_0}(G)$$

      We conclude that if $0<\epsilon_1< d/3$, and if $G$ is a $\epsilon_1-$perturbation of $F$, then the atoms $\widehat A \in {\cal A}_{k_0} (G)$  remain at  distance larger than $d/3>0$ from the  separation line $\widehat S$ of $G$, and at Hausdorff distance smaller than $\epsilon _1 <d/3$ of its corresponding atom $A \in {\cal A}_{k_0} (F)$, being $d = \dist (K_0, S)$ and
 $ K_0 = \cup \{A \in {\cal A}_{k_0}(F)\}  = F^{k_0}(B)  $.

      Recall that $$k \geq k_0 \;\; \Rightarrow \; \; K = F^k(B) \subset F^{k_0}(B) = K_0 , \;\; \; G^k(B) \subset G^{k_0}(B)$$ $$ A '\in {\cal A}_k (F)  \Rightarrow   A '\subset  A \in {\cal A}_{k_0}(F),\; \;  \widehat A '\in {\cal A}_k (G)  \Rightarrow  \widehat A '\subset \widehat A \in {\cal A}_{k_0}(G)$$ $$ \Rightarrow  \dist (\widehat A ', \widehat S) \geq \dist (\widehat A, \widehat S) \geq d = \dist (K_0, S). $$ So, using the same positive numbers $d >0$ and $0 <\epsilon_1 < d/3$ for all $k \geq k_0$, we deduce from the hypothesis of Lemma \ref{lema}, the following statement:

      {\bf (S) } There exists $k_0 \geq 1$ and $0 <\epsilon_1 < d/3$,  such that if $G$ is a $\epsilon_1-$perturbation of $F$, then for all $k \geq k_0$ the atoms $\widehat A \in {\cal A}_{k} (G)$  remain at  distance larger than $d /3>0$ from the  separation line $\widehat S$ of $G$, and at Haussdorf distance smaller than $\epsilon _1 < d /3$ of its corresponding atom $A \in {\cal A}_{k} (F)$, being $d = \dist (K_0 , S), \; \; K_0 = \cup \{A \in {\cal A}_{k_0 }(F)\}  = F^{k_0}(B) $.

    Now,  $k_0 \geq 1$, $\epsilon_1  >0$ and $d >0$ are fixed as in statement (S), and the generation $k \geq k_0$ is chosen and also \em fixed, \em  such that the atoms of ${\cal A}_k(F)$ and of ${\cal A}_k(G)$ have all diameter smaller than $d/6$.

    Repeating the argument (P) used in the proof of Lemma \ref{lema}, we deduce that $A \in {\cal A}_k(F)$ is  in the interior of some (and unique) continuity piece $B_{i_k}$ of $F$:
     $A \subset int (B_{i_k})$, and
    all the points at distance smaller than $d/3$ of $A$   are contained in $int  (B_{i_k}) $. This last includes the  atom $\widehat A = \Psi (A) \in {\cal A}_{k}(G)$. Then  $\widehat A \subset int (B_{i_k}).$  Repeating once more the same argument (P) used in the proof of Lemma \ref{lema}, now with $G$ instead of $F$, we deduce  $\widehat A \subset int (\widehat B_{\widehat i_k})$
      for some unique $\widehat i_k$.

      {\bf (T)} We assert  that for the fixed $k \geq k_0$ constructed as above, for any $A \in {\cal A}_k(F), \; \; \widehat A = \Psi(A)  \in {\cal A}_k(G)$, the indexes $i_k$ and $\widehat i_k$ constructed as above, coincide: $\widehat i_k = i_k$.

        By contradiction, if  $\widehat i_k \neq i_k$ then, applying Proposition \ref{remarkagregado},  the distance from any point $p \in B_{i_k} \cap B_{\widehat i_k}$ to  $\widehat S$ is smaller than $\mbox{ Hdist }(S, \widehat S) < \epsilon_1 $. Then $\mbox{dist}(A, \widehat S) <  \epsilon _1 < d/3$ contradicting the statement (S).  $\; \; \; \; \Box$ {\bf{(T)}}

 So $i_k$ is the first index of the itinerary of the atom $A \in {\cal A}_k (F)$, which due to (T)  coincides with the first index $\widehat i_k$ of the itinerary of the atom $\widehat A = \Psi (A) \in {\cal A}_k(G)$. Now let us prove that  the indexes of the  itinerary of $A$ and of $\widehat A$, i.e. the indexes for their future iterates, also coincide.

    The future iterates of any atom of generation $k$, is an atom of generation $k' \geq k$ by $F$, and also by $G$. They are contained in some atoms of generation $k$ of $F$, and $G$ respectively. Therefore, using (T), the  images  of an atom $A \in {\cal A}_k(F)$ or $\widehat A = \Psi (A) \in {\cal A}_k (G)$,  by  all the future iterates of $F$ or of $G$ respectively, are in the interior of their respective one-to-one corresponding continuity pieces $B_{i_{k'}} , \; \widehat B_{i_{k'}}  $, where the index $i_{k'}$ is the same for all $k' \geq k$. Then the itineraries of $A$ and $\widehat A = \Psi (A)$ are the same, as we asserted in { (Q)}. $\; \; \; \Box$ {\bf (Q)}

    As a consequence of assertion (Q),  the indexes $i_0, i_1, i_2, \ldots$ in the finite chain of atoms denoted in (\ref{equationImagenesDeAtomo}) and (\ref{ecuacionChain2}),  remain unchanged, for $F$ or for $G$, being $G$ an $\epsilon-$perturbation of $F$ for $\epsilon >0$ small enough.  We deduce the following statement:

     \textbf{A}: \em  The number of periodic orbits in the atoms of generation $k$,  and their periods,  remain unchanged, when substituting $F$ by any $\epsilon$-perturbation $G$, if $\epsilon >0$ is sufficiently small. \em

    Now it is standard to prove by induction on $k \geq 1$ the following property:

     Let $F$ be a piecewise continuous contractive map with contraction rate upper bounded by $0 <\lambda <1$. Let $\epsilon >0$ such that $\lambda + \epsilon = \widehat \lambda  < 1$. Let $G$  be an $\epsilon-$ perturbation of $F$. Then, for all $k \geq 1$, each atom $\widehat A$ of generation $k$ for $G$, is at distance smaller than $\sum_{j= 0}^{k-1} 2 \epsilon \, \widehat \lambda  ^j < 2 \epsilon/ (1 - \widehat \lambda)  = \epsilon ^* >0$ of the respective atom $A$ of generation $k$ for $F$,  with the same itinerary than $\widehat A$.

    Therefore we deduce the following statement:

     \textbf{B}: \em Any periodic point found in an atom $\widehat A  $ of generation $k$ for $G$, is at distance smaller than $\epsilon ^*$ than the respective periodic point found in the corresponding atom $A$ for $F$ with the same itinerary.  \em

    The statements A and B imply that the limit cycles are
    persistent according to Definition \ref{definicionPersistencia}.$\;\; \Box$

 \begin{remark} \em
In the proof of Lemmas \ref{lema} and \ref{lemma2}, we did not use the separation property $f_i(B_i) \cap f_j(B_j) = \emptyset\; \; \forall i \neq j $. At the very beginning of the proof of Lemma \ref{lemma2}, we obtained that  the piecewise continuous and locally contractive systems verifying the thesis of the Lemma \ref{lema}, \ even if they do not have the separation property, contain an open family of systems  in the topology defined in \ref{topologia}. Then:

\em In the space of all the piecewise continuous and locally contractive systems (even if they do not have the separation property), those whose limit set is formed by a finite number of persistent limit cycles form  an open family. \em

 Nevertheless, to prove the genericity of the periodic persistent behavior, we  need to prove that the family of periodic maps is dense in the space of systems. In the  proof  of the Theorem \ref{Teorema1}, to obtain the density property, we shall restrict to the space of systems ${\cal S}$   that  verify the separation property.

 \end{remark}

 \begin{remark}
\label{remarkperiodosgrandes} \em
From the proof of Lemma \ref{lema},  the first integer $k_0 \geq 1$ such that $F^{k_0} (B) \bigcap S = \emptyset$ may be as large as wanted. Therefore, there exists periodic systems such that the number of iterates that takes to settle into a periodic behavior, from any initial state, may be arbitrarily large. In particular, if the initial state is that of a periodic orbit, the minimum number of iterates that takes to return to it, i.e. the period $p$, may be arbitrarily large.

This fact has a relevant consequence in the   applications to experimental Science: For instance, if a network of $n \approx 10^{12}$ neurons, is such that  no  neuron becomes dead, i.e. it does not eventually remain forever under the threshold level without completing its oscillatory cycle, then the periodic sequences $ i_1, \ldots, i_p$, defined in the proof of Lemma \ref{lemma2}, as the itinerary of the periodic limit cycles, during their period $p$, have at least once each of all the indexes $i \in \{1, 2, \ldots, n\}$. Then $p \geq n \approx 10^{12}$.

Besides, as it is proved in \cite{Yo}, there exists a minimum time $T >0$ between two consequent spikes (instantaneous pulses coupling the neurons of the network), if the neuronal network is completely inhibitory. Accordingly to the values of the electric components of each neuron, it is reasonably to assume for instance that $T \approx 10 \; [ms]$. As $n \approx 10^{12}$, the lasting time of the periodic sequence  could be approximately $ 10^{-3} \times 10^{12} [s] = 10^9 [s] \geq 31$ years. So, if most of the neurons of an artificial network were inhibitory\footnote{\lq\lq From the neuroscience point of view, the model of a population of inhibitory (biological) coupled pacemaker neurons is very hypothetical." Personal communication from a reader of the preliminary version of this paper.}, and if most of them did not become dead, the observation of the theoretical periodic behavior of the inhibitory system in the future, could not be practical during a reasonable time of experimentation, if the electronic circuits (and thus the time constants) are designed with values in a scale 1/1 respect to the electronic model of a biological neuron. Therefore only the irregularities inside the period could be registered by the observer, showing the system as virtually chaotic (\cite{cessac}). In other words, the rate of the amount of information of the network during a   relatively short time of experimentation, will be positive, even if the mathematical entropy, computed as the limit of the rate of information when the time  goes to infinite, is  zero in the generic periodic dynamics.
\end{remark}

{\bf \em Proof of Theorem \ref{Teorema1}. } Due to Lemma \ref{lema} the existence of a finite number of limit cycles attracting all the orbits of the space is verified at least for those systems in the hypothesis of \ref{lema}. This hypothesis is an open condition because $K_0 = F^{k_0}(B)$ and $S$ are compact set at positive distance, and for fix $k_0$, the set $F^{k_0}(B)$ depends continuously on the map $F$.

To prove its genericity  it is enough to prove  that the hypothesis of Lemma \ref{lema} is also a dense condition in the space ${\cal S}$ of piecewise continuous contractive maps with the separation property, with the topology in ${\cal S}$ defined in \ref{topologia}.

Take $F$ being  not
  finally periodic.

  We shall prove that, for all $\epsilon >0$ there exists a $\epsilon-$ perturbation $G$ of $F$ that verifies the hypothesis of Lemma \ref{lema}, and thus $G$ is finally periodic with persistent limit cycles.

Let be given an arbitrarily  small $\epsilon >0$.

The contractive homeomorphisms $f_i$ of the finite family
  $F = \{f_i: B_i \mapsto B\}_i $, with contraction rate $0 < \lambda <1, $
  can be $C^0$  extended to $$F_{\displaystyle{\epsilon}}= \{f_{i, \epsilon}: U_i \mapsto B\}_i, \; \ \; \mbox{ where }\; \ \; f_{i, \epsilon}: U_i \mapsto B, \; \; \; f_{i, \epsilon}|_{B_i} = f_i, $$ $U_i$ is a compact neighborhood  such that $ B_i = \overline B_i \subset int (U_i) \subset U_i = \overline U_i \subset B,$ and $f_{i, \epsilon}$ is an homeomorphism onto its image.

   We construct  $f_{i, \epsilon}$  still contractive in $U_i$,  with a  contraction rate \begin{equation}
  \label{equationLoli0}
  0 < \lambda ' < 1 \mbox{ such that } |\lambda -\lambda'| < \epsilon.\end{equation}

  Such a finite family $F_{\displaystyle {\epsilon}} $ of continuous extensions $f_{i, \epsilon}$ to open sets $U_i \supset B_i$, exists as an application of Tietze Theorem (see for instance Theorem 2.15 of \cite{armstrong}), applied to homeomorphisms.

The role of the family $F_{\displaystyle {\epsilon}} $ of continuous extensions $f_{i, \epsilon}$
will be the following:

The union of the domains of $f_{i, \epsilon}$ is the union of the  sets $U_i \subset B$. They do not form  a partition of $B$ because they overlap on  sets with non void interiors, covering the discontinuity line $S$ of the given $F$. So $F_{\displaystyle{\epsilon}}$ is multi-defined now, not only in $S$ but in the   set $$V=\bigcup_{i \neq j} U_i \cap U_j \supset S $$ with non void interior.  The covering $\{U_i\}$ makes the line of discontinuities $S$ a kind of fuzzy set: i.e. one can move freely the line of discontinuities $S $ inside the interior of the  set $V$, to define a new partition of the space $B$.

Our purpose is to find some $G$ that is a $\epsilon$- perturbation of $F$, such that $G$  verifies the hypothesis of Lemma \ref{lema}. We will choose not any $G$, but one in a very particular way, obtained from $F$ moving \em only \em the line $S$ of discontinuities of $F$ to a new line $S_{\widehat{{\cal P}}} \subset V$, and the partition ${\cal P}$ of continuity pieces of $F$ to a near new partition ${\widehat{{\cal P}}}$. We will do that \em without changing the functional values of $F$ in the points where it was already defined. \em

The image of $B$ by the future    $n$-th. iterate of $F_{\displaystyle{\epsilon}}$, includes the image of $B$ by $F^n$, because $f_{i, \epsilon}$ is defined in a  set $U_i \supset B_i$ (recall that $B_i$ is the domain of $f_i$), and   $f_{i, \epsilon}|_{B_i} = f_i$.

 But the image of $B$ by the future $n -$th. iterate  of $F_{\displaystyle{\epsilon}}$, includes also the image of $B$ by   $G^n$  (being $G$ any piecewise contractive function $G$ that is a restriction of $F_{\displaystyle{\epsilon}}$ to some continuity pieces $C_i \subset U_i$).
 Then, $F_{\displaystyle{\epsilon}}^n(B)$ includes the image of $B$ by the iterate of all those $\epsilon $- perturbation $G$ of $F$, obtained from $F$ moving only its line of discontinuities, and so, changing only the partition ${\cal P} = \{B_i\}$ of the continuity pieces  to a new partition ${\widehat{{\cal P}}}=\{C_i\}$  such that $C_i \subset U_i$   (without changing the functional values of $F$ in $B_i \cap C_i$).

 In other words, the extended family $F_{\displaystyle{\epsilon}}$ is the ``egg" of all the $\epsilon-$ perturbations $G$ of $F$, obtained from $F$ moving \em only \em the partition ${\cal P}$ to a new partition ${\widehat{{\cal P}}}$, that is, moving the line of discontinuities $S$ to a new line $S_{\widehat{{\cal P}}} $ (contained in the   set where $F_{\displaystyle{\epsilon}}$ is multidefined).

 The extended
 map $F_{\displaystyle{\epsilon} } = \{f_{i, \epsilon}: U_i \mapsto B\}_i $, is now multidefined in $\bigcup _{i \neq j } U_i \cap U_j \supset S$.
  The separation property is an open condition, thus the extension $F_{\displaystyle{\epsilon}}$ still verifies $f_{i, \epsilon} (U_i ) \cap f_{j, \epsilon} (U_j) = \emptyset $ for all $i \neq j $,  if the neighborhoods $U_i$ and $U_j$ are chosen at a sufficiently small Hausdorff distance from their respective pieces $B_i$ and $B_j$, and  $\epsilon>0$ is small enough.

 Call $\epsilon_1>0$ to a positive real number  smaller or equal than $\epsilon$, and also smaller or equal than
  the distance from
 $B_i$ to the complement of $U_i$, for all
 $i = 1,2, \ldots m$. Precisely \begin{equation}
 \label{equationLoli}
 0 < \epsilon_1   = \min \{\epsilon, \; \; \min_{1 \leq i \leq m} \dist (B_i, U_i^c)  \}\end{equation}

 Consider the compact sets:
 \begin{equation}
\label{ecuacionConjuntoK}
K = \bigcap_{k \geq 1}\; \; \;  \bigcup _{(i_1, \ldots, i_k) \in \{1,2 \ldots m\}^k}
f_{i_k} \circ \ldots \circ f_{i_1}(B_{i_1})\end{equation}
$$K^+ = \bigcap_{k \geq 1}\; \; \;  \bigcup _{(i_1, \ldots, i_k) \in \{1,2 \ldots m\}^k}
f_{i_k, \epsilon} \circ \ldots \circ f_{i_1, \epsilon}(U_{i_1}) \;\;\; \supset \;\;
K$$
The sets $K$ and $K^+$ are the forward limit sets of $F$ and $F_{\epsilon}$ respectively.

Define the family ${\cal A}_{k, \epsilon}\;  $ \em  of the extended atoms \em of generation $k\geq 1$ for $F_{\displaystyle{\epsilon}} $ that form $K^+$, defined as follows:

The set $A \subset B$ is an extended atom of generation $k \geq 1$ if and only if there exists a word $(i_1, \ldots,  i_{k-1} ,  i_k) \in \{1, 2, \ldots, m\} ^k$ such that $$A= f_{i_k, \epsilon} \circ \ldots \circ f_{i_1, \epsilon}(U_{i_1}).$$
The diameter of each extended atom of generation $k$ is
 smaller that $diam (B) \cdot  {(\lambda')} ^k$ because $f_{i, \epsilon}$ is contractive with contraction rate $0 < \lambda ' < 1$. Therefore, for sufficiently large $k\geq 1$ all
the extended atoms of generation $k$ that form $K^+$ have diameters  smaller that $\epsilon_1 /2 $:
\begin{equation}
\label{equationLoli2} A \in {\cal A}_{k, \epsilon } \; \; \Rightarrow \; \; \mbox{diam}(A) < \frac{\epsilon_1}{2}.
\end{equation}

We assert that the extended atoms of  generation $k \geq 1$ are pairwise disjoint: in fact,  for two different $i \neq j $ the images  are disjoint: $f_{i, \epsilon}(U_i) \cap f_{j, \epsilon} (U_j)= \emptyset$. So the atoms of generation 1 are pairwise disjoint. Two extended atoms of generation $k \geq1$ are $f_{i_k, \epsilon} \circ \ldots \circ f_{i_1, \epsilon}(U_{i_1})$ and $f_{j_k, \epsilon} \circ \ldots \circ f_{j_1, \epsilon}(U_{j_1})$. They can intersect if and only if $(i_1, i_2, \ldots, i_k) =(j_1, j_2, \ldots, j_k)$ because each $f_{i, \epsilon}$  is an homeomorphism onto its image. So, they intersect if and only if they coincide.

By construction, $U_i \supset B_i$ and $f_{i, \epsilon}|_{B_i} = f_i$. Therefore each of the  atoms of generation $k$ for $F$, is contained in the respective extended atom of generation $k$ for $F_{\displaystyle{\epsilon}}$, that has the same finite word $(i_1, i_2, \ldots, i_k) $.

If none of the extended  atoms of generation $k$ intersects $S$, then none of the atoms of generation $k$ for $F$ intersects $S$, and the system verifies the hypothesis of Lemma \ref{lema}. So, in this case, there is nothing to prove, because $F$ is finally periodic. (Recall our assumption at the beginning of this proof that the given $F$ is not finally periodic.)

On the other hand, if some of
the extended atoms of generation $k$ intersects $S$, consider a new finite partition ${\widehat{{\cal P}}}= \{C_i\}_{1 \leq i \leq m}$ of $B$ such that the distance, defined in (\ref{definicionHausdorfDistance}),
between ${\widehat{{\cal P}}}$ and the given partition ${\cal P }$ of $F$, is smaller
than $\epsilon_1 >0$: \begin{equation}
\label{equationLoli3}
 \dist ({\cal P}, \; {\widehat{{\cal P}}} ) < \epsilon_1 \leq \epsilon, \end{equation}
 where $\epsilon _1 >0$ was defined in the equality (\ref{equationLoli}).

Choose the new partition ${\widehat{{\cal P}}}$ such that the new separation line  $S_{\widehat{{\cal P}}} = \bigcup _{i
\neq j}
 (C_i \cap C_j) $ does not intersect the extended atoms of generation $k$ of $K^+$:
   \begin{equation}
   \label{equationLoli4}
   S_{\widehat{{\cal P}}} \bigcap \left( \bigcup \{A \in {\cal A}_{k, \epsilon}\}\right ) \; = \; \emptyset \end{equation}

   This last condition is possible because the diameters of the extended atoms $A \in {\cal A}_{k, \epsilon}$ are all smaller than $\epsilon_1 /2$, due to inequality (\ref{equationLoli2}). They are compact pairwise disjoint sets, because of the separation property. The distance between the two partitions ${\cal P}$ and ${\widehat{{\cal P}}}$  is smaller than $\epsilon_1>0$ due to inequality (\ref{equationLoli3}), but  can be chosen larger than $\epsilon_1 /2$, and such that does not cut the atoms $A \in {\cal A}_{k, \epsilon}$, which verify inequality (\ref{equationLoli2}) and are all pairwise disjoint compact sets.

   Due to the construction above and to the definition in equality (\ref{definicionHausdorfDistance}), the maximum  Hausdorff distance between the  respective pieces $B_i$ of ${\cal P}$ and $C_i$ of ${\widehat{{\cal P}}}$ is larger than $\epsilon_1/2 >0$ and smaller than $\epsilon_1>0$.

   We note that the old, and principally the new, separation lines $S_{{\cal P}}$ and $S_{{\widehat{{\cal P}}}}$,  are not necessarily  $C^1$ nor even Lipschitz manifolds in the space $B$, and even if they are, they do not need to be $\epsilon _1$- $C^1$ or Lipschitz near one from the other, to be $\epsilon_1$ near with the Hausdorff distance.

  The
 condition $\dist ({\cal P}, {\widehat{{\cal P}}}) < \epsilon_1$ in (\ref{equationLoli3}), joined with the assumption

 \noindent $\dist (U_i^c, B_i) \geq \epsilon_1 $ in (\ref{equationLoli}), where $B_i$ is the $i-$th piece of the partition ${\cal P}$, implies that the respective piece $C_i$ of the partition ${\widehat{{\cal P}}}$ verifies $$C_i \subset U_i.$$ Therefore the
extension $f_{i, \epsilon}: U_i \mapsto B$ in $ F_{\displaystyle{\epsilon}}$ whose domain of definition is $U_i$  can be restricted to $C_i$.

 Define $$G  = \{g_i: C_i \mapsto B\}_{1 \leq  m} \mbox{ where } g_i =
 f_{i, \epsilon}|_{C_i}.$$ By construction $G$ and $F$ coincide in $C_i \cap B_i$, the distance between the respective partitions ${\cal P}$ and ${\widehat{{\cal P}}}$ is smaller than $\epsilon_1 \leq \epsilon$ due to (\ref{equationLoli3}), and the difference of their respective contraction rates $\lambda '$ and $\lambda$ is also smaller than $\epsilon$, due to (\ref{equationLoli0}). So $G$ is a $\epsilon  $-perturbation
 of the given $F$, according to the Definition  \ref{topologia}.

 It is enough now, to prove  that $G$ is finally periodic with persistent limit cycles.

 Consider the limit set $K_G$ of $G$ as follows:
 $$K_{G} = \bigcap_{k \geq 1} \; \; \;  \bigcup _{(i_1, \ldots, i_k) \in \{1,2 \ldots m\}^k}
g_{i_k} \circ \ldots \circ g_{i_1}(C_{i_1})$$

As $G$ is a restriction of $F_{\displaystyle{\epsilon} } $ to the sets $C_i \subset U_i$, we
have  that $K_G \subset K^+$, and  in particular for all $k \geq 1$ the atoms of generation $k$ for $G$, i.e. $g_{i_k} \circ \ldots \circ g_{i_1}(C_{i_1})$, are contained in the extended atoms of generation $k$ for $F_{\displaystyle{\epsilon}}$.

Due to inequality (\ref{equationLoli4}), the separation line $S_{G} = S_{{\widehat{{\cal P}}}}$ among the continuity pieces $C_i$ of $G$ is disjoint with the extended atoms of generation $k$ of $F_{\displaystyle{\epsilon}}$. Therefore, it is also disjoint with the atoms  of generation $k$ of $G$. Then  $G^k(B) \bigcap S_G =
\emptyset$ and, applying lemma \ref{lema}, $G$ is finally
periodic with persistent limit cycles. $\Box$

It is possible (but not immediate) to construct, in a  compact ball $B$ of any dimension $n \geq 2$, piecewise continuous systems, uniformly locally contractive and with the separation property, as defined in Section \ref{definicionesabstractas}, that do not verify the thesis of Lemma  \ref{lema}, and thus their limit set $K$, defined by the Equality (\ref{ecuacionConjuntoK}) is not composed only by periodic limit cycles, but contains a Cantor set attractor.

%The first author was partially financed by Project PDT54/001.

%

\end{document}